\documentclass[notitlepage,leqno,10pt]{amsart}
\usepackage{amsmath,amssymb,amsbsy,amsfonts,amsthm,latexsym,
            amsopn,amstext,amsxtra,euscript,amscd}
\usepackage{hyperref}

\begin{document}
\bibliographystyle{plain}
\newfont{\teneufm}{eufm10}
\newfont{\seveneufm}{eufm7}
\newfont{\fiveeufm}{eufm5}
%
%
\newfam\eufmfam
              \textfont\eufmfam=\teneufm \scriptfont\eufmfam=\seveneufm
              \scriptscriptfont\eufmfam=\fiveeufm
\def\bbbr{{\rm I\!R}}
\def\bbbm{{\rm I\!M}}
\def\bbbn{{\rm I\!N}}
\def\bbbf{{\rm I\!F}}
\def\bbbh{{\rm I\!H}}
\def\bbbk{{\rm I\!K}}
\def\bbbp{{\rm I\!P}}
\def\bbbone{{\mathchoice {\rm 1\mskip-4mu l} {\rm 1\mskip-4mu l}
{\rm 1\mskip-4.5mu l} {\rm 1\mskip-5mu l}}}
\def\bbbc{{\mathchoice {\setbox0=\hbox{$\displaystyle\rm C$}\hbox{\hbox
to0pt{\kern0.4\wd0\vrule height0.9\ht0\hss}\box0}}
{\setbox0=\hbox{$\textstyle\rm C$}\hbox{\hbox
to0pt{\kern0.4\wd0\vrule height0.9\ht0\hss}\box0}}
{\setbox0=\hbox{$\scriptstyle\rm C$}\hbox{\hbox
to0pt{\kern0.4\wd0\vrule height0.9\ht0\hss}\box0}}
{\setbox0=\hbox{$\scriptscriptstyle\rm C$}\hbox{\hbox
to0pt{\kern0.4\wd0\vrule height0.9\ht0\hss}\box0}}}}
\def\bbbq{{\mathchoice {\setbox0=\hbox{$\displaystyle\rm
Q$}\hbox{\raise
0.15\ht0\hbox to0pt{\kern0.4\wd0\vrule height0.8\ht0\hss}\box0}}
{\setbox0=\hbox{$\textstyle\rm Q$}\hbox{\raise
0.15\ht0\hbox to0pt{\kern0.4\wd0\vrule height0.8\ht0\hss}\box0}}
{\setbox0=\hbox{$\scriptstyle\rm Q$}\hbox{\raise
0.15\ht0\hbox to0pt{\kern0.4\wd0\vrule height0.7\ht0\hss}\box0}}
{\setbox0=\hbox{$\scriptscriptstyle\rm Q$}\hbox{\raise
0.15\ht0\hbox to0pt{\kern0.4\wd0\vrule height0.7\ht0\hss}\box0}}}}
\def\bbbt{{\mathchoice {\setbox0=\hbox{$\displaystyle\rm.
T$}\hbox{\hbox to0pt{\kern0.3\wd0\vrule height0.9\ht0\hss}\box0}}
{\setbox0=\hbox{$\textstyle\rm T$}\hbox{\hbox
to0pt{\kern0.3\wd0\vrule height0.9\ht0\hss}\box0}}
{\setbox0=\hbox{$\scriptstyle\rm T$}\hbox{\hbox
to0pt{\kern0.3\wd0\vrule height0.9\ht0\hss}\box0}}
{\setbox0=\hbox{$\scriptscriptstyle\rm T$}\hbox{\hbox
to0pt{\kern0.3\wd0\vrule height0.9\ht0\hss}\box0}}}}
\def\bbbs{{\mathchoice
{\setbox0=\hbox{$\displaystyle     \rm S$}\hbox{\raise0.5\ht0\hbox
to0pt{\kern0.35\wd0\vrule height0.45\ht0\hss}\hbox
to0pt{\kern0.55\wd0\vrule height0.5\ht0\hss}\box0}}
{\setbox0=\hbox{$\textstyle        \rm S$}\hbox{\raise0.5\ht0\hbox
to0pt{\kern0.35\wd0\vrule height0.45\ht0\hss}\hbox
to0pt{\kern0.55\wd0\vrule height0.5\ht0\hss}\box0}}
{\setbox0=\hbox{$\scriptstyle      \rm S$}\hbox{\raise0.5\ht0\hbox
to0pt{\kern0.35\wd0\vrule height0.45\ht0\hss}\raise0.05\ht0\hbox
to0pt{\kern0.5\wd0\vrule height0.45\ht0\hss}\box0}}
{\setbox0=\hbox{$\scriptscriptstyle\rm S$}\hbox{\raise0.5\ht0\hbox
to0pt{\kern0.4\wd0\vrule height0.45\ht0\hss}\raise0.05\ht0\hbox
to0pt{\kern0.55\wd0\vrule height0.45\ht0\hss}\box0}}}}
\def\bbbz{{\mathchoice {\hbox{$\sf\textstyle Z\kern-0.4em Z$}}
{\hbox{$\sf\textstyle Z\kern-0.4em Z$}}
{\hbox{$\sf\scriptstyle Z\kern-0.3em Z$}}
{\hbox{$\sf\scriptscriptstyle Z\kern-0.2em Z$}}}}
\def\ts{\thinspace}

\newtheorem{theorem}{Theorem}
\newtheorem{lemma}[theorem]{Lemma}
\newtheorem{claim}[theorem]{Claim}
\newtheorem{cor}[theorem]{Corollary}
\newtheorem{prop}[theorem]{Proposition}
\newtheorem{definition}[theorem]{Definition}
\newtheorem{remark}[theorem]{Remark}
\newtheorem{question}[theorem]{Open Question}
\newtheorem{example}[theorem]{Example}

\def\qed{\ifmmode
\squareforqed\else{\unskip\nobreak\hfil
\penalty50\hskip1em\null\nobreak\hfil\squareforqed
\parfillskip=0pt\finalhyphendemerits=0\endgraf}\fi}

\def\squareforqed{\hbox{\rlap{$\sqcap$}$\sqcup$}}

\def \C {{\mathbb C}}
\def \F {{\mathbb F}}
\def \L {{\mathbb L}}
\def \K {{\mathbb K}}
\def \Q {{\mathbb Q}}
\def \Z {{\mathbb Z}}
\def\cA{{\mathcal A}}
\def\cB{{\mathcal B}}
\def\cC{{\mathcal C}}
\def\cD{{\mathcal D}}
\def\cE{{\mathcal E}}
\def\cF{{\mathcal F}}
\def\cG{{\mathcal G}}
\def\cH{{\mathcal H}}
\def\cI{{\mathcal I}}
\def\cJ{{\mathcal J}}
\def\cK{{\mathcal K}}
\def\cL{{\mathcal L}}
\def\cM{{\mathcal M}}
\def\cN{{\mathcal N}}
\def\cO{{\mathcal O}}
\def\cP{{\mathcal P}}
\def\cQ{{\mathcal Q}}
\def\cR{{\mathcal R}}
\def\cS{{\mathcal S}}
\def\cT{{\mathcal T}}
\def\cU{{\mathcal U}}
\def\cV{{\mathcal V}}
\def\cW{{\mathcal W}}
\def\cX{{\mathcal X}}
\def\cY{{\mathcal Y}}
\def\cZ{{\mathcal Z}}
\newcommand{\rmod}[1]{\: \mbox{mod}\: #1}

\def\tcN{\cN^\mathbf{c}}
\def\F{\mathbb F}
\def\Tr{\operatorname{Tr}}
\def\mand{\qquad \mbox{and} \qquad}
\renewcommand{\vec}[1]{\mathbf{#1}}
\def\eqref#1{(\ref{#1})}
\newcommand{\ignore}[1]{}
\hyphenation{re-pub-lished}
\parskip 1.5 mm
\def\lln{{\mathrm Lnln}}
\def\Res{\mathrm{Res}\,}
\def\F{{\bbbf}}
\def\Fp{\F_p}
\def\fp{\Fp^*}
\def\Fq{\F_q}
\def\ff{\F_2}
\def\ffn{\F_{2^n}}
\def\K{{\bbbk}}
\def \Z{{\bbbz}}
\def \N{{\bbbn}}
\def\Q{{\bbbq}}
\def \R{{\bbbr}}
\def \P{{\bbbp}}
\def\Zm{\Z_m}
\def \Um{{\mathcal U}_m}
\def \Bf{\frak B}
\def\Km{\cK_\mu}
\def\va {{\mathbf a}}
\def \vb {{\mathbf b}}
\def \vc {{\mathbf c}}
\def\vx{{\mathbf x}}
\def \vr {{\mathbf r}}
\def \vv {{\mathbf v}}
\def\vu{{\mathbf u}}
\def \vw{{\mathbf w}}
\def \vz {{\mathbfz}}
\def\\{\cr}
\def\({\left(}
\def\){\right)}
\def\fl#1{\left\lfloor#1\right\rfloor}
\def\rf#1{\left\lceil#1\right\rceil}
\def\flq#1{{\left\lfloor#1\right\rfloor}_q}
\def\flp#1{{\left\lfloor#1\right\rfloor}_p}
\def\flm#1{{\left\lfloor#1\right\rfloor}_m}
\def\Al{{\sl Alice}}
\def\Bob{{\sl Bob}}
\def\Or{{\mathcal O}}
\def\inv#1{\mbox{\rm{inv}}\,#1}
\def\invM#1{\mbox{\rm{inv}}_M\,#1}
\def\invp#1{\mbox{\rm{inv}}_p\,#1}
\def\Ln#1{\mbox{\rm{Ln}}\,#1}
\def \nd {\,|\hspace{-1.2mm}/\,}
\def\ord{\mu}
\def\E{\mathbf{E}}
\def\Cl{{\mathrm {Cl}}}
\def\epp{\mbox{\bf{e}}_{p-1}}
\def\ep{\mbox{\bf{e}}_p}
\def\eq{\mbox{\bf{e}}_q}
\def\bm{\bf{m}}
\newcommand{\floor}[1]{\lfloor {#1} \rfloor}
\newcommand{\comm}[1]{\marginpar{
\vskip-\baselineskip
\raggedright\footnotesize
\itshape\hrule\smallskip#1\par\smallskip\hrule}}
\def\rem{{\mathrm{\,rem\,}}}
\def\dist {{\mathrm{\,dist\,}}}
\def\etal{{\it et al.}}
\def\ie{{\it i.e. }}
\def\veps{{\varepsilon}}
\def\eps{{\eta}}
\def\ind#1{{\mathrm {ind}}\,#1}
               \def \MSB{{\mathrm{MSB}}}
\newcommand{\abs}[1]{\left| #1 \right|}

\title{On a problem of Pillai with Fibonacci numbers and powers of $2$}
\author{
{\sc Mahadi Ddamulira, Florian Luca and  Mihaja Rakotomalala}
}
\address{
AIMS Ghana (Biriwa)\\ P.O. Box DL 676\\ Adisadel, Cape Coast\\ Central Region, Ghana\\
}
\email{mahadi@aims.edu.gh\\
mihaja@aims.edu.gh} 
\address{
School of Mathematics\\ 
University of the Witwatersrand \\
Private Bag 3, Wits 2050, Johannesburg, South Africa \\
}
\email{Florian.Luca@wits.ac.za}

\begin{abstract}
In this paper, we  find all all integers $c$ having at least two representations as a difference between a Fibonacci number and a power of $2$.
\end{abstract} 

\maketitle

\section{Introduction} 

Let $\{F_n\}_{n\ge 0} $ be the sequence of Fibonacci numbers given by $F_0=0$, $F_1=1$ and 
$$
F_{n+2}=F_{n+1}+F_n\quad  {\text{\rm for all}}\quad n\geq 0.
$$
Its  first few terms are
$$
1,1,2,3,5,8,13,21,34,55,89,144,233,377,610,987,1597,2584,4181,\ldots
$$
In 1936 (see \cite{Pi1, Pi2}), Pillai showed that if $a$ and $b$ are coprime then there exists $c_0(a,b)$ such that if $c>c_0(a,b)$ is an integer, then the equation $c=a^x-b^y$ 
has at most one positive integer solution $(x,y)$. In the special case $(a,b)=(3,2)$ which was studied before Pillai by Herschfeld \cite{Her1, Her2}, Pillai conjectured that the only integers $c$ admitting 
two representations of the form $3^x-2^y$ are given by
$$
1=3-2=3^2-2^3,\qquad -5=3-2^3=3^3-2^5,\qquad -13=3-2^4=3^5-2^8.
$$
This was confirmed by R. J. Stroeker and Tijdeman in 1982 (see \cite{ST}). Here we study a related problem and find all positive integers $c$ admitting two representations of the form $F_n-2^m$ 
for some positive integers $n$ and $m$. We assume that representations with $n\in \{1,2\}$ (for which $F_1=F_2$) count as one representation just to avoid trivial ``parametric families" such as
$1-2^m=F_1-2^m=F_2-2^m$, and so we always assume that $n\ge 2$. Notice the solutions
\begin{eqnarray}
\label{eq:sol}
1& = & 5-4=3-2(=F_5-2^2=F_4-2^1),\nonumber\\
-1 & = & 3-4=1-2(=F_4-2^2=F_2-2^1),\nonumber\\
-3 &= & 5-8=1-4=13-16(=F_5-2^3=F_2-2^2=F_7-2^4),\nonumber\\
5 & = & 21-16=13-8(=F_8-2^4=F_7-2^3),\nonumber\\
0 & = & 8-8=2-2(=F_6-2^3=F_3-2^1),\\
-11 & = & 21-32=5-16(=F_8-2^5=F_5-2^4),\nonumber\\
-30 & = & 34-64=2-32(=F_9-2^6=F_3-2^5)\nonumber\\
85 & = & 4181-4096 = 89-4(=F_{19}-2^{12}=F_{11}-2^2).\nonumber
\end{eqnarray}
We prove the following theorem.

\begin{theorem}
\label{thm:1}
The only integers $c$ having at least two representations of the form $F_n-2^m$ are $c\in \{0,1,-1,-3,5,-11,-30,85\}$. Furthermore, for each $c$ in the above set, all its representations of the form $F_n-2^m$ with 
integers $n\ge 2$ and $m\ge 1$ appear in the list \eqref{eq:sol}.
\end{theorem}

\section{A lower bound for a linear forms in logarithms of algebraic numbers}
\label{sec:2}

In this section, we state a result concerning lower bounds for linear
forms in logarithms of algebraic numbers, which will be used in the
proof of our theorem.

Let $\eta$ be an algebraic number of degree $d$, whose
minimal polynomial over the integers is $$g(x) = a_0 \prod_{i=1}^d (x - \eta^{(i)}).$$
The logarithmic height of $\eta$ is defined as
$$
h(\eta) = \frac{1}{d}\( \log |a_0| + \sum_{i=1}^d \log \max \{ |\eta^{(i)}|,
1\}\).
$$
Let $\L$ be an algebraic number field and $d_{\L}$ be the degree of
the field $\L$. Let $\eta_1, \eta_2, \ldots, \eta_l \in \L$  not $0$
or $1$ and $d_1, \ldots, d_l$ be nonzero integers. We put
$$
D =\max\{|d_1|, \ldots, |d_l|, 3\},
$$
and put
$$
\Lambda = \prod_{i=1}^l \eta_i^{d_i} -1.
$$
Let $A_1, \ldots, A_l$ be positive integers such that $$A_j \geq h'(\eta_j) := \max \{d_{\L}h(\eta_j), |\log \eta_j|, 0.16
\}\quad {\text{\rm for}}\quad j=1,\ldots, l.
$$
The following result is due to Matveev \cite{Matveev}.
\begin{theorem}
\label{thm:Matveev} If $\Lambda \neq 0$ and $\L \subset \R $, then
\begin{equation*}
\label{ineq:matveev} \log |\Lambda| > -1.4 \cdot
30^{l+3}l^{4.5}d_{\L}^2(1+ \log d_{\L})(1+ \log D)A_1A_2\cdots A_l.
\end{equation*}
\end{theorem}

\section{\small Proof of Theorem \ref{thm:1}}

Assume that $(n,m)\ne (n_1,m_1)$ are such that 
$$
F_n-2^m=F_{n_1}-2^{m_1}.
$$
If $m=m_1$, then $F_n=F_{n_1}$ and since $\min\{n,n_1\}\ge 2$, we get that $n=n_1=2$, so $(n,m)=(n_1,m_1)$, which is not the case. Thus, $m\ne m_1$, and we may assume that $m>m_1$. Since 
\begin{equation}
\label{eq:main}
F_n-F_{n_1}=2^m-2^{m_1},
\end{equation}
and the right--hand side is positive, we get that the left--hand side is also positive and so $n>n_1$. Thus, $n\ge 3$ and $n_1\ge 2$. We use the Binet formula
$$
F_k=\frac{\alpha^k-\beta^k}{\alpha-\beta}\qquad {\text{\rm for~all}}\qquad k\ge 0,
$$
where $(\alpha,\beta)=((1+{\sqrt{5}})/2,(1-{\sqrt{5}})/2)$ are the roots of the characteristic equation $x^2-x-1=0$ of the Fibonacci sequence. It is well-known that
$$
\alpha^{k-2}\le F_k\le \alpha^{k-1}\qquad {\text{\rm for~all}}\qquad k\ge 1.
$$
In \eqref{eq:main} we have
\begin{eqnarray}
\label{eq:t}
\alpha^{n-4} & \le & F_{n-2}\le F_n-F_{n_1}=2^m-2^{m_1}<2^m,\\
\alpha^{n-1} & \ge & F_n>F_n-F_{n_1}=2^m-2^{m_1}\ge 2^{m-1}\nonumber,
\end{eqnarray}
therefore
\begin{equation}
\label{eq:c1}
n-4<c_1 m\quad {\text{\rm and}}\quad n-1>c_1 (m-1),\quad {\text{\rm where}}\quad c_1=\log 2/\log \alpha=1.4402\ldots.
\end{equation}
If $n<400$, then $m<300$. We ran a computer program for $2\le n_1<n\le 400$ and $1\le m_1<m<300$ and found only the solutions from list \eqref{eq:sol}. From now, on, $n\ge 400$.
By the above inequality \eqref{eq:c1}, we get that $n>m$. 
Thus, we get
\begin{eqnarray*}
\left|\frac{\alpha^n}{\sqrt{5}}-2^m\right| & = & \left|\frac{\beta^{n}}{\sqrt{5}}+\frac{\alpha^{n_1}-\beta^{n_1}}{\sqrt{5}}-2^{m_1}\right|\le \frac{\alpha^{n_1}+2}{\sqrt{5}}+2^{m_1}\\
& \le &  
\frac{2\alpha^{n_1}}{\sqrt{5}}+2^{m_1}<2\max\{\alpha^{n_1},2^{m_1}\}.
\end{eqnarray*}
Dividing by $2^m$ we get
\begin{equation}
\label{eq:lin}
\left| {\sqrt{5}}^{-1} \alpha^n 2^{-m}-1\right|<2\max\left\{\frac{\alpha^{n_1}}{2^m}, 2^{m_1-m}\right\}<\max\{\alpha^{n_1-n+6} ,2^{m_1-m+1}\},
\end{equation}
where for the last right--most inequality above we used \eqref{eq:t} and the fact that $2<\alpha^2$.
For the left--hand side above, we use Theorem \ref{thm:Matveev} with the data
$$
l=3,~\eta_1={\sqrt{5}},~\eta_2=\alpha,~\eta_3=2,~d_1=-1,~d_2=n,~d_3=-m.
$$
We take ${\mathbb L}={\mathbb Q}({\sqrt{5}})$ for which $d_{\mathbb L}=2$. Then we can take $A_1=2h(\eta_1)=\log 5,$ $A_2=2h(\eta_2)=\log \alpha,~A_3=2h(\eta_3)=2\log 2$.  We take $D=n$. We have
$$
\Lambda={\sqrt{5}}^{-1} \alpha^n 2^{-m}-1.
$$
Clearly, $\Lambda\ne 0$, for if $\Lambda=0$, then $\alpha^{2n}\in {\mathbb Q}$, which is false. The left--hand side of \eqref{eq:lin1} is bounded, by Theorem \ref{thm:Matveev}, as 
$$
\log |\Lambda|>-1.4\times 30^6 \times 3^{4.5} \times 2^2 (1+\log 2) (1+\log n) (\log 5) (2\log \alpha) (2\log 2).
$$
Comparing with \eqref{eq:lin}, we get
$$
\min\{(n-n_1-6)\log \alpha,(m-m_1-1)\log 2\}<1.1\times 10^{12} (1+\log n),
$$
which gives
$$
\min\{(n-n_1)\log \alpha, (m-m_1)\log 2\}<1.2\times 10^{12} (1+\log n).
$$
Now the argument splits into two cases.

\medskip

{\bf Case 1.} $\min\{(n-n_1)\log \alpha, (m-m_1)\log 2\}=(n-n_1)\log \alpha$.

\medskip

In this case, we rewrite \eqref{eq:main} as
$$
\left|\frac{(\alpha^{n-n_1}-1)}{\sqrt{5}} \alpha^{n_1} -2^m\right|=\left|\frac{\beta^n-\beta^{n_1}}{\sqrt{5}}-2^{m_1}\right|< 2^{m_1}+1\le 2^{m_1+1},
$$
so 
\begin{equation}
\label{eq:lin1}
\left|\left(\frac{\alpha^{n-n_1}-1}{\sqrt{5}}\right) \alpha^{n_1} 2^{-m}-1\right|<2^{m_1-m-1}.
\end{equation}

\medskip

{\bf Case 2.} $\min\{(n-n_1)\log \alpha, (m-m_1)\log 2\}=(m-m_1)\log 2$.

\medskip

In this case, we rewrite \eqref{eq:main} as
$$
\left|\frac{\alpha^n}{\sqrt{5}}-2^{m_1}(2^{m-m_1}-1)\right|=\left|\frac{\beta^n+\alpha^{n_1}-\beta^{n_1}}{\sqrt{5}}\right|< \frac{\alpha^{n_1}+2}{\sqrt{5}}<\alpha^{n_1},
$$
so 
\begin{equation}
\label{eq:lin2}
\left|({\sqrt{5}} (2^{m-m_1}-1))^{-1} \alpha^{n} 2^{-m_1}-1\right|<\frac{\alpha^{n_1}}{2^{m}-2^{m_1}}\le \frac{2\alpha^{n_1}}{2^m}\le 2\alpha^{n_1-n+4}<\alpha^{n_1-n+6}.
\end{equation}

\medskip

Inequalities \eqref{eq:lin1} and \eqref{eq:lin2} suggest studying lower bounds for the absolute values of
$$
\Lambda_1=\left(\frac{\alpha^{n-n_1}-1}{\sqrt{5}}\right) \alpha^{n_1} 2^{-m}-1\quad {\text{\rm and}}\quad \Lambda_2=({\sqrt{5}} (2^{m-m_1}-1))^{-1} \alpha^{n} 2^{-m_1}-1.
$$
We apply again Theorem \ref{thm:Matveev}. We take in both cases $l=3,~Ê\eta_2=\alpha,~\eta_3=2$. For $\Lambda_1$, we have $d_2=n_1,~d_3=-m$, while for $\Lambda_2$ we have $d_2=n,~d_3=-m_1$. In both cases
we take $D=n$. We take 
$$
\eta_1=\frac{\alpha^{n-n_1}-1}{\sqrt{5}},\quad {\text{\rm or}}\quad \eta_1={\sqrt{5}} (2^{m-m_1}-1),
$$
according to whether we work with $\Lambda_1$ or $\Lambda_2$, respectively. For $\Lambda_1$ we have $d_1=1$ and for $\Lambda_2$ we have $d_1=-1$. In both cases ${\mathbb L}={\mathbb Q}({\sqrt{5}})$ for which $d_{\mathbb L}=2$. The minimal polynomial of $\eta_1$ divides 
$$
5X^2-5F_{n-n_1}X-((-1)^{n-n_1}+1-L_{n-n_1})\quad {\text{\rm or}}\quad X^2-5(2^{m-m_1}-1)^2,
$$
respectively, where $\{L_k\}_{k\ge 0}$ is the Lucas companion sequence of the Fibonacci sequence given by $L_0=2,~L_1=1$, $L_{k+2}=L_{k+1}+L_k$ for $k\ge 0$ for which 
its Binet formula of the general term is 
$$
L_k=\alpha^k+\beta^k\qquad {\text{\rm for ~all}}\qquad k\ge 0.
$$
Thus,
\begin{equation}
\label{eq:uuu}
h(\eta_1)\le \frac{1}{2} \left(\log 5+\log\left(\frac{\alpha^{n-n_1}+1}{\sqrt{5}}\right)\right)\quad {\text{\rm or}}\quad \log({\sqrt{5}}(2^{m-m_1}-1),
\end{equation}
respectively. In the first case,
\begin{equation}
\label{eq:vvv}
h(\eta_1)<\frac{1}{2} \log(2{\sqrt{5}} \alpha^{n-n_1})<\frac{1}{2} (n-n_1+4)\log \alpha<7\times 10^{11} (1+\log n),
\end{equation}
and in the second case 
$$
h(\eta_1)<\log(8\times 2^{m-m_1})=(m-m_1+3)\log 2<1.3\times 10^{12} (1+\log n).
$$
So, in both cases, we can take $A_1=2.6\times 10^{12} (1+\log n).$ We have to justify that $\Lambda_i\ne 0$ for $i=1,2$. But $\Lambda_1=0$ means
$$
(\alpha^{n-n_1}-1)\alpha^{n_1}={\sqrt{5}}\times  2^m.
$$
Conjugating this relation in ${\mathbb Q}$, we get that
\begin{equation}
\label{eq:alphabeta}
(\alpha^{n-n_1}-1)\alpha^{n_1}=-(\beta^{n-n_1}-1)\beta^{n_1}.
\end{equation}
The absolute value of the left-hand side is at least $\alpha^{n}-\alpha^{n_1}\ge \alpha^{n-2}\ge \alpha^{398}$, while the absolute value of the right--hand side is at most 
$(|\beta|^{n-n_1}+1)|\beta|^{n_1}<2$, which is a contradiction. As for $\Lambda_2$, note that $\Lambda_2=0$ implies $\alpha^{2n}\in {\mathbb Q}$, which is not possible. 
 We then get that 
$$
\log |\Lambda_i|>-1.4\times 30^6\times 3^{4.5} \times 2^2 (1+\log 2) (1+\log n) (2.6\times 10^{12} (1+\log n)) 2(\log 2)\log\alpha,
$$ 
for $i=1,2$. Thus,
$$
\log |\Lambda_i|>-1.7\times 10^{24} (1+\log n)^2\quad {\text{\rm for}}\quad i=1,2.
$$
Comparing these with \eqref{eq:lin1} and \eqref{eq:lin2}, we get that
$$
(m-m_1-1)\log 2<1.7 \times 10^{24} (1+\log n)^2,\quad (n-n_1-6)\log \alpha<1.7\times 10^{24} (1+\log n)^2,
$$
according to whether we are in Case 1 or in Case 2. Thus, in both Case 1 and Case 2, we have 
\begin{eqnarray}
\label{eq:minmax}
\min\{(n-n_1)\log \alpha, (m-m_1)\log 2\} & < & 1.2\times 10^{12} (1+\log n)\nonumber\\ 
\max\{(n-n_1)\log \alpha,(m-m_1)\log 2\} & < & 1.8 \times 10^{24} (1+\log n)^2.
\end{eqnarray}
We now finally rewrite equation \eqref{eq:main} as 
$$
\left|\frac{(\alpha^{n-n_1}-1)}{\sqrt{5}} \alpha^{n_1} -2^{m_1}(2^{m-m_1}-1)\right|=\left|\frac{\beta^n-\beta^{n_1}}{\sqrt{5}}\right|<|\beta|^{n_1}=\frac{1}{\alpha^{n_1}}.
$$
We divide both sides above by $2^{m}-2^{m_1}$ getting
\begin{eqnarray}
\label{eq:lin3}
\left|\left(\frac{\alpha^{n-n_1}-1}{{\sqrt{5}}(2^{m-m_1}-1)}\right) \alpha^{n_1} 2^{-m_1}-1\right| & < & \frac{1}{\alpha^{n_1} (2^m-2^{m_1})}\le \frac{2}{\alpha^{n_1} 2^m}\nonumber\\
& \le &  
2\alpha^{4-n-n_1}\le \alpha^{4-n},
\end{eqnarray}
because $\alpha^{n_1}\ge \alpha^2>2$. To find a lower-bound on the left--hand side above, we use again Theorem \ref{thm:Matveev} with the data
$$
l=3,~\eta_1=\frac{\alpha^{n-n_1}-1}{{\sqrt{5}}(2^{m-m_1}-1)},~\eta_2=\alpha,~\eta_3=2,~d_1=1,~d_2=n_1,~d_3=-m_1,~D=n.
$$
We have ${\mathbb L}={\mathbb Q}({\sqrt{5}})$ with $d_{\mathbb L}=2$. Using that $h(x/y)\le h(x)+h(y)$ for any two nonzero algebraic numbers $x$ and $y$, we have
\begin{eqnarray*}
h(\eta_1) & \le & h\left(\frac{\alpha^{n-n_1}-1}{\sqrt{5}}\right)+h(2^{m-m_1}+1)< \log(2{\sqrt{5}} \alpha^{n-n_1})+\log( 2^{m-m_1}+1)\\
& \le &  
(n-n_1)\log \alpha+(m-m_1)\log 2+\log (2{\sqrt{5}})+1<2\times 10^{24} (1+\log n)^2,
\end{eqnarray*}
where in the above chain of inequalities we used the arguments from \eqref{eq:uuu} and \eqref{eq:vvv} as well as the bound \eqref{eq:minmax}.
So, we can take $A_1=4\times 10^{24} (1+\log n)^2$ and certainly $A_2=\log \alpha$ and $A_3=2\log 2$. We need to show that if we put 
$$
\Lambda_3= \frac{(\alpha^{n-n_1}-1)}{\sqrt{5}} \alpha^{n_1} -2^{m_1}(2^{m-m_1}-1),
$$
then $\Lambda_3\ne 0$. But $\Lambda_3=0$ leads to 
$$
(\alpha^{n-n_1}-1)\alpha^{n_1}={\sqrt{5}} (2^m-2^{m_1}),
$$
which upon conjugation in ${\mathbb L}$ leads to \eqref{eq:alphabeta}, which we have seen that it is impossible. Thus, $\Lambda_3\ne 0$. Theorem \ref{thm:Matveev} gives
$$
\log |\Lambda_3|>-1.4\times 30^6 \times 3^{4.5} \times 2^2 (1+\log 2) (1+\log n) (4\times 10^{24} (1+\log n)^2) 2(\log 2)\log \alpha,
$$
which together with \eqref{eq:lin3}  gives
$$
(n-3)\log \alpha<3\times 10^{36} (1+\log n)^3,
$$
leading to $n<7\times 10^{42}$.

Now we need to reduce the bound. To do so, we make use several times of the following result, which is a slight variation of a result due to Dujella and Peth\H{o} which itself is a generalization of a result of Baker and Davenport \cite{BD}. For a real number $x$, we put  $ ||x|| = \min\{|x-n|\,:\, n \in \Z\}$ for the distance from $x$ to the nearest integer.
\begin{lemma}
\label{Dujella-Petho}
Let $M$ be a positive integer, let $p/q$ be a convergent of the continued
fraction of the irrational
$\tau$ such that $q > 6M$, and let $A, B,  \mu$ be some real
numbers with $A > 0$ and $B > 1$. Let $\varepsilon :=||\mu q|| - M||\tau q||$. If $\varepsilon > 0$, then there is no solution to the
inequality
$$
0 < m\tau - n + \mu < AB^{-k},
$$
in positive integers $m, n$ and $k$ with
$$
m \leq M \qquad {\text{and}} \qquad k \geq \dfrac{\log(Aq/\varepsilon)}{\log B}.
$$
\end{lemma}

We return first to \eqref{eq:lin} and put 
$$
\Gamma=n\log \alpha-m\log 2-\log {\sqrt{5}}.
$$
Assume that $\min\{n-n_1,m-m_1\}\ge 20$ and we go to \eqref{eq:lin}. This is not a very restrictive assumption since, as we shall see immediately, if this condition fails then we do the following:
\begin{itemize}
\item[(i)] if $n-n_1<20$ but $m-m_1\ge 20$, we go to \eqref{eq:lin1};
\item[(ii)] if $n-n_1\ge 20$ but $m-m_1<20$, we go to \eqref{eq:lin2};
\item[(iii)] if both $n-n_1<20$ and $m-m_1<20$, we go to \eqref{eq:lin3}.
\end{itemize}

In \eqref{eq:lin}, since $|e^{\Gamma}-1|=|\Lambda|<1/4$, we get that $|\Gamma|<1/2$. Since $|x|<2|e^x-1|$ holds for all $x\in (-1/2,1/2)$, we get that 
$$
|\Gamma|<2\max\{\alpha^{n_1-n+6},2^{m-m_1+1}\}\le \max\{\alpha^{n_1-n+8}, 2^{m_1-m+2}\}.
$$
Assume $\Gamma>0$. Then
$$
0<n\left(\frac{\log \alpha}{\log 2}\right)-m+\frac{\log(1/{\sqrt{5}})}{\log 2}<\max\left\{\frac{\alpha^{8}}{(\log 2) \alpha^{n-n_1}}, \frac{4}{(\log 2) 2^{m-m_1}}\right\}.
$$
We apply Lemma \ref{Dujella-Petho} with 
$$
\tau=\frac{\log \alpha}{\log 2},\quad \mu=\frac{\log(1/{\sqrt{5}})}{\log 2},\quad (A,B)=(68,\alpha)\quad {\text{\rm or}}\quad (6,2).
$$
We let $\tau=[a_0,a_1,a_2,\ldots]=[0,1,2,3,1,2,3,2,4,\ldots]$ be the continued fraction of $\tau$. We take $M=7\times 10^{42}$. We take 
$$
\frac{p}{q}=\frac{p_{149}}{q_{149}}={\frac{\scriptstyle 75583009274523299909961213530369339183941874844471761873846700783141852920}{\scriptstyle108871285052861946543251595260369738218462010383323482629611084407107090003}}
$$
where $q>10^{74}>6 M$. We have $\varepsilon>0.09$, therefore either
$$
n-n_1\le \frac{\log(68q/0.09)}{\log\alpha}<369\quad {\text{\rm or}}\quad m-m_1\frac{\log(6q/0.09)}{\log 2}<253.
$$
Thus we have that either $n-n_1\le 368$ or $m-m_1\le 252$. A similar conclusion is obtained when $\Gamma<0$. 

In case $n-n_1\le 368$, we go to \eqref{eq:lin1}. There, we assume that $m-m_1\ge 20$. We put 
$$
\Gamma_1=n_1\log \alpha -m\log 2+\log\left(\frac{\alpha^{n-n_1}-1}{\sqrt{5}}\right).
$$
Then \eqref{eq:lin1} implies that 
$$
|\Gamma_1|<\frac{4}{2^{m-m_1}}.
$$
Assume $\Gamma_1>0$. Then
$$
0<n_1\left(\frac{\log \alpha}{\log 2}\right)-m+\frac{\log((\alpha^{n-n_1}-1)/{\sqrt{5}})}{\log 2}<\frac{4}{(\log 2) 2^{m-m_1}}<\frac{6}{2^{m-m_1}}.
$$
We keep the same $\tau,~M,~q$, $(A,B)=(6,2)$ and put
$$
\mu_k=\frac{\log((\alpha^k-1)/{\sqrt{5}})}{\log 2},\qquad k=1,2,\ldots, 368.
$$
We have problems at $k\in \{4,12\}$. We discard these values and we will treat them later. For the remaining values of $k$, we get $\varepsilon>0.001$. Hence, 
by Lemma \ref{Dujella-Petho}, we get
$$
m-m_1<\frac{\log(6 q/0.001)}{\log 2}<259.
$$
Thus, $n-n_1\le 368$ implies $m-m_1\le 258$, unless $n-n_1\in \{4,12\}$.  A similar conclusion is reached if $\Gamma_1<0$ with the same two exceptions for $n-n_1\in \{4,12\}$. The reason we have a problem at $k\in \{4,8\}$ is because
$$
\frac{\alpha^4-1}{\sqrt{5}}=\alpha^2\quad {\text{\rm and}}\quad \frac{\alpha^{12}-1}{\sqrt{5}}= 2^3\alpha^6.
$$
So,
$$
\Gamma_1=(n_1+2) \tau-m,\quad {\text{\rm or}}\quad (n_1+6)\tau-(m-3)\quad {\text{\rm when}}\quad k=4,12,\quad {\text{\rm respectively}}.
$$
Thus we get that
$$
\left|\tau-\frac{m}{n_1+2}\right|<\frac{4}{2^{m-m_1} (n_1+2)}\quad {\text{\rm or}}\quad \left|\tau-\frac{m-3}{n_1+6}\right|<\frac{4}{2^{m-m_1} (n_1+6)}.
$$
Assume $m-m_1>150$. Then $2^{m-m_1}> 8\times (8\times 10^{42})> 8\times (n_1+6)$, therefore 
$$
\frac{4}{2^{m-m_1} (n_1+2)}<\frac{1}{2(n_1+2)^2}\quad {\text{\rm and}}\quad \frac{4}{2^{m-m_1}(n_1+6)}<\frac{1}{2(n_1+6)^2}.
$$
By a criterion of Legendre, it follows that $m/(n_1+2)$ or $(m+3)/(n_1+6)$ are convergents of $\tau$, respectively. So, say one of $m/(n_1+2)$ or $m/(n_1+6)$ is of the form $p_k/q_k$ for some $k=0,1,2,\ldots,99$. Here we use that $q_{99}>8\times 10^{42}>n_1+6$. Then 
$$
\frac{1}{(a_k+2) q_k^2}<\left|\tau-\frac{p_k}{q_k}\right|.
$$
Since $\max\{a_k:k=0,\ldots,99\}=134$, we get that 
$$
\frac{1}{136 q_k^2}< \frac{4}{2^{m-m_1} q_k}\quad {\text{\rm and}}\quad q_k\quad {\text{\rm divides~one~of}}\quad \{n_1+2,~n_1+6\}.
$$
Thus
$$
2^{m-m_1}\le 4\times 136 (n_1+6)<4\times 136\times 8\times 10^{42}
$$
giving $m-m_1\le 151$. Hence, even in the case $n-n_1\in \{4,12\}$, we still keep the conclusion that $m-m_1\le 258$. 

Now let us assume that $m-m_1\le 252$. Then we go to \eqref{eq:lin2}. We write 
$$
\Gamma_2=n\log \alpha -m_1\log 2+\log(1/({\sqrt{5}}(2^{m-m_1}-1))).
$$
We assume that $n-n_1\ge 20$. Then 
$$
|\Gamma_2|<\frac{2\alpha^6}{\alpha^{n-n_1}}.
$$
Assuming $\Gamma_2>0$, we get
$$
0<n\left(\frac{\log \alpha}{\log 2}\right)-m_1+\frac{\log(1/({\sqrt{5}}(2^{m-m_1}-1)))}{\log 2}<\frac{2\alpha^6}{(\log 2) \alpha^{n-n_1}}<\frac{52}{\alpha^{n-n_1}}.
$$
We apply again Lemma \ref{Dujella-Petho} with the same $\tau,~q,~M$, $(A,B)=(52,\alpha)$ and 
$$
\mu_k=\frac{\log(1/({\sqrt{5}}(2^k-1))}{\log 2}\quad {\text{\rm for}}\quad k=1,2,\ldots,252.
$$ 
We get $\varepsilon>0.0005$, therefore 
$$
n-n_1<\frac{\log(52 q/0.0005)}{\log \alpha}< 379.
$$
A similar conclusion is reached when $\Gamma_2<0$. To conclude, we first got that either $n-n_1\le 368$ or $m-m_1\le 252$. If $n-n_1\le 368$, then $m-m_1\le 258$, and if $m-m_1\le 252$, then $n-n_1\le 378$. In conclusion, we always have $n-n_1<380$ and $m-m_1<260$. 

Finally we go to \eqref{eq:lin3}. We put
$$
\Gamma_3=n_1\log \alpha-m_1\log 2+\log\left(\frac{\alpha^{n-n_1}-1}{{\sqrt{5}}(2^{m-m_1}-1)}\right).
$$
Since $n\ge 400$, \eqref{eq:lin3} tells us that 
$$
|\Gamma|<\frac{2}{\alpha^{n-3}}=\frac{2\alpha^3}{\alpha^n}.
$$
Assume that $\Gamma_3>0$. Then
$$
0<n_1 \left(\frac{\log \alpha}{\log 2}\right)-m_1+\frac{\log((\alpha^k-1)/{\sqrt{5}}(2^\ell-1))}{\log 2}<\frac{2\alpha^3}{(\log 2) \alpha^n}<\frac{13}{\alpha^n}
$$
where $(k,l):=(n-n_1,m-m_1)$. We apply again Lemma \ref{Dujella-Petho} with the same $\tau,~M,~q$, $(A,B)=(13,\alpha)$ and 
$$
\mu_{k,l}=\frac{\log((\alpha^k-1)/{\sqrt{5}}(2^\ell-1))}{\log 2}\quad {\text{\rm for}}\quad 1\le k\le 379,~1\le l\le 259. 
$$
We have a problem at $(k,l)=(4,1),~(12,1)$ (as for the case of \eqref{eq:lin1}) and additionally for $(k,l)=(8,2)$ since
$$
\frac{\alpha^8-1}{{\sqrt{5}}(2^2-1)}=\alpha^4.
$$
We discard the cases $(k,l)=(4,1),~(12,1),(8,2)$ for the time being. For the remaining ones, we get $\varepsilon>7\times 10^{-6}$, so we get
$$
n\le \frac{\log(13 q/(7\times 10^{-6}))}{\log \alpha}<385.
$$
A similar conclusion is reached when $\Gamma_3<0$. Hence $n<400$. Now we look at the cases $(k,l)=(4,1),~(12,1),~(8,2)$. The cases $(k,l)=(4,1),~(12,1)$ can be treated as we did before 
when we showed that $n-n_1\le 368$ implies  $m-m_1\le 258$. The case when $(k,l)=(8,2)$ can be dealt with similarly as well. Namely, it gives
$$
|(n_1+4)\tau-m_1|<\frac{13}{\alpha^n}.
$$
Hence
\begin{equation}
\label{eq:www}
\left|\tau-\frac{m_1}{n_1+4}\right|<\frac{13}{(n_1+4)\alpha^n}.
\end{equation}
Since $n\ge 400$, then $\alpha^n>2\times 13\times (8\times 10^{42})>2\times 13 (n_1+4)$, which shows that the right--hand side of inequality \eqref{eq:www} is at most $2/(n_1+4)^2$. By Legendre's criterion, $m/(n_1+4)=p_k/q_k$ for some $k=0,1,\ldots,99$. We then get by an argument similar to a previous one that
$$
\alpha^n\le 13\times 136 \times (8\times 10^{42})
$$
giving $n\le 220$. So, the conclusion is that $n<400$ holds also in the case of the pair $(k,l)=(8,2)$. However, this contradicts our working assumption that $n\ge 400$. 

Theorem \ref{thm:1} is therefore proved.  \qed

\section*{Acknowledgements} We thank the referee for comments which improved the quality of this paper. We also thank J. J. Bravo for spotting a computational oversight in a previous version of this paper.


\begin{thebibliography}{99}

\bibitem{BD} A. Baker and H. Davenport,``The equations $3x^2-2=y^2$ and $8x^2-7=z^2$", {\it Quart.J of Math. Ser.(2)\/} {\bf 20} (1969), 129--137.


\bibitem{DP} A.~Dujella and A.~Peth\H o, ``A generalization of a theorem of
Baker and Davenport", {\it Quart. J. Math. Oxford Ser. (2)\/} {\bf
49} (1998), 291--306.

\bibitem{Her1} A. Herschfeld, ``The equation $2^x-3^y=d$", {\it Bull. Amer. Math. Soc.\/} {\bf 41} (1935), 631.

\bibitem{Her2}  A. Herschfeld, ``The equation $2^x-3^y=d$", {\it Bull. Amer. Math. Soc.\/} {\bf 42} (1936), 231--234.


\bibitem{Matveev}E. M. Matveev, ``An explicit lower bound for a homogeneous rational linear form in the logarithms of algebraic numbers II", {\it Izv. Ross. Akad. Nauk Ser. Mat.\/}  {\bf 64} (2000), 125--180; translation in Izv. Math. {\bf 64} (2000), 1217--1269.


\bibitem{Pi1} S. S. Pillai,  `` On $a^x - b^y = c$", {\it J. Indian math. Soc. (N.S.)\/}  {\bf 2} (1936), 119--122.

\bibitem{Pi2} S. S. Pillai, ``A correction to the paper ÒOn $a^x - b^y = c$Ó", {\it J. Indian math. Soc.\/} {\bf 2} (1937),
215.

\bibitem{ST} R. J. Stroeker and R. Tijdeman, ``Diophantine equations", in {\it Computational methods
in number theory\/}, Part II, vol. {\bf 155} of Math. Centre Tracts, Math. Centrum, Amsterdam,
1982, 321Ð369.

\end{thebibliography}
\end{document}